\documentclass[final]{siamltex}
\usepackage{amsfonts}

\if@onethmnum
  \newtheorem{example}[theorem]{Example}
\else
  \newtheorem{example}[theorem]{Example}
\fi

\title{Operational Method for Finite Difference Equations}

\author{Salvador Merino C\'ordoba\thanks{Applied Mathematics Department,  University 
        of M\'alaga, M\'alaga, Spain. smerino@uma.es}}
 
\begin{document}

\maketitle

\begin{abstract}
In this article I present a fast and direct method for solving several types of linear finite difference 
equations (FDE) with constant coefficients. The method is based on a polynomial form of the translation 
operator and its inverse, and can be used to find the particular solution of the FDE. This work raises the possibility of developing new ways to expand the scope of the operational methods.
\end{abstract}



\pagestyle{myheadings}
\thispagestyle{plain}
\markboth{Salvador Merino C\'ordoba}{OPERATIONAL METHOD FOR FINITE DIFFERENCE EQUATIONS}

\section[Finite Difference Equations]{Finite Difference Equations (FDEs)}

These equations describe the relationship between the present
value of a function and a discrete set of $n$ previous values:
\[f(t),f(t+1),\dots,f(t+n)\]

The solution of an FDE is a function $f(t)$, with $t\in
\mathbb{Z}$, that satisfies the equation. 

Typically, the $n$ known values for a FDE of degree $n$ are referred to as 
initial conditions. In order to establish a unique solution, it is
necessary to know the values of the initial conditions as precisely as possible. 
They are habitually defined at equal intervals starting from
$t=0$: $f(0),f(1),f(2),\dots,f(n)$. In this notation, the FDE predicts
$f(n)$.

\subsection{Linear FDEs} A linear FDE can be expressed as follows:
\[a_0(t)y(t+n)+a_{1}(t)y(t+n-1)+\cdots+a_{n-1}(t)y(t+1)+a_n(t)y(t)=\phi(t)\]
To solve a {\em linear equation with constant
coefficients} (LECC),
\[a_0y(t+n)+a_{1}y(t+n-1)+\cdots+a_{n-1}y(t+1)+a_n y(t)=\phi(t) \ ,\]
we define the \textbf{Translation Operator} \index{operator!translation} $T$:  
    \[ y(t+1)=Ty(t); \quad y(t+k)=\underbrace{T(T(\cdots(T}_{k}y(t))\cdots))=T^k y(t) \Rightarrow \] \vspace{-.6cm}
    \[ \Rightarrow a_0 T^n y(t)+a_1 T^{n-1}y(t)+\cdots+a_{n-1}Ty(t)+a_n y(t)=\phi(t) \Rightarrow P(T)y(t)=\phi(t)\]
and look for a general solution of the form
\begin{equation}
y_{General}(t) = y_{Homogenous}(t)+y_{Particular}(t) \Rightarrow y_G(t)=y_H(t)+y_P(t) \ .
    \label{EDF:solutiongeneral}
\end{equation}

Several methods for solving LECCs exist. In this paper, I develop an operational method 
to find a particular solution.

\section[Operational Method for FDE]{Operational Method for Finite Difference Equations}

The basic idea of the operational method is simple. If we wish to
define the polynomial operator $P(T)$, is is possible to
establish its properties? We begin with the following axiom:

\begin{equation} \mathbf{
P(T)g(t)=f(t) \Leftrightarrow \frac{1}{P(T)}f(t)=g(t)}
\label{EDF:axioma}
\end{equation}

\subsection{Properties of $\displaystyle{\frac{1}{P(T)}}$}
\begin{enumerate}
    \item \textbf{Linearity}
    \[P(T)\Bigl[\alpha f(t)+\beta g(t)\Bigl]=\]
    \[=a_0\Bigl[\alpha f(t+n)+\beta g(t+n)\Bigl]+a_1\Bigl[\alpha f(t+n-1)+\beta g(t+n-1)\Bigl]+\]
    \[+\cdots+a_{n-1}\Bigl[\alpha f(t+1)+\beta g(t+1)\Bigl]+a_n\Bigl[\alpha f(t)+\beta g(t)\Bigl]=\]
    \[=\alpha\Bigl[a_0f(t+n)+a_1 f(t+n-1)+\cdots+a_{n-1}f(t+1)+a_n f(t)\Bigl]+\]
    \[+\beta\Bigl[a_0 g(t+n)+a_1 g(t+n-1)+\cdots+a_{n-1}g(t+1)+a_n g(t)\Bigl]=\]
    \[=\alpha P(T) f(t)+\beta P(T) g(t)\]
    \\
    Therefore:
    \[P(T)\Bigl[\alpha f(t)+\beta g(t)\Bigl]=\alpha P(T) f(t)+\beta P(T) g(t)\Leftrightarrow \]
     \begin{equation}
         \Leftrightarrow \mathbf{ \frac{1}{P(T)}\Bigl[\alpha f(t)+\beta g(t)\Bigl]=
         \alpha \frac{1}{P(T)} f(t)+\beta \frac{1}{P(T)} g(t)}
         \label{EDF:Linealidad}
      \end{equation}

    \item \textbf{Inverse Translation}
    \[ T^n f(t)=f(t+n) \Rightarrow f(t)=\frac{1}{T^n}f(t+n) \Rightarrow\]
    We change $t \to t-n$: 
    \begin{equation}
    \displaystyle{\Rightarrow\mathbf{
    \frac{1}{T^n}f(t)=f(t-n)}}
    \label{EDF:TranslacionInversa}
    \end{equation}

    \item \textbf{Unity}
    \\
    Let $y(t)=t$, with $n\in \mathbb{N}$. The expression $y(t+1)-y(t)$
    is equivalent to $t+1-t=1$, so the FDE $y(t+1)-y(t)=1$
    has the solution $y(t)=t$. Finally,
     \begin{equation}
    \displaystyle\mathbf{{\frac{1}{T-1}(1)=t}}
    \label{EDF:Unidad}
    \end{equation}

    \item \textbf{Propagation}
    \\
    Let us solve the two simplest cases directly:
    \[y(t+1)-y(t)=t \Rightarrow \frac{t}{T-1}=\frac{t}{1}\cdot\frac{t-1}{2}\]
    \[y(t+2)-2y(t+1)+y(t)=t \Rightarrow \frac{t}{(T-1)^2}=\frac{t}{1}\cdot\frac{t-1}{2}\cdot\frac{t-2}{3} \ .\]
    Generalizing the progression, we obtain
    \[\frac{t}{(T-1)^n}=\frac{t}{1}\cdot\frac{t-1}{2}\cdot\frac{t-2}{3}\cdots\frac{t-n}{n+1}=\prod_{i=0}^n \frac{t-i}{i+1} \ .\]
    Therefore, 
    \begin{equation}
    \displaystyle{\mathbf{\frac{t}{(T-1)^n}=\prod_{i=0}^n \frac{t-i}{i+1}}}
    \label{EDF:Propagacion}
    \end{equation}

\end{enumerate}

\subsection{Equations where $\mathbf{\phi(t)=\lambda^t}$}
Consider the polynomial 
\[P(T)=a_0 T^n+a_1T^{n-1}+\cdots+a_{n-1}T+a_n \ .\]

We have
\[P(T)\lambda^t=a_0
T^n\lambda^t+a_1T^{n-1}\lambda^t+\cdots+a_{n-1}T\lambda^t+a_n\lambda^t=\]
\[=a_0
\lambda^{t+n}+a_1\lambda^{t+n-1}+\cdots+a_{n-1}\lambda^{t+1}+a_n\lambda^t=\]
\[=\lambda^t(a_0
\lambda^{n}+a_1\lambda^{n-1}+\cdots+a_{n-1}\lambda+a_n)=\lambda^tP(\lambda)\]
\\
Therefore, given that $P(\lambda) \ne 0$ always holds, we consider the
property
\begin{equation}
       P(T)\lambda^t = \lambda^tP(\lambda) \Leftrightarrow \mathbf{\frac{1}{P(T)}\lambda^t=
       \frac{\lambda^t}{P(\lambda)}}
        \label{EDF:Potencia}
\end{equation}

\begin{example}
Find a particular solution $y_P(t)$ of the finite difference
equation
\[y(t+2)-5y(t+1)+4y(t)=3^t \ . \]
\end{example}
\begin{solution}
\[y(t+2)-5y(t+1)+4y(t)=3^t\Rightarrow (T^2-5T+4)y(t)=3^t\]
Replacing:
\[ y(t)=\frac{1}{T^2-5T+4}3^t=\frac{3^t}{3^2-5 \cdot 3+4} \Rightarrow y_P(t)=-\frac{3^t}{2}\]
\end{solution}

\subsection{Equations where $\mathbf{\phi(t)=\cos n \pi t}$ or $\mathbf{\phi(t)=\sin n \pi t}$}

We begin with Euler's formula:
\[ e^{n \pi t i}=\cos n \pi t + i \sin n \pi t \ . \]
In particular, we consider $n\in\mathbb{N}$:
\[ e^{n \pi i}=\cos n \pi + i \sin n \pi = (-1)^n + i \cdot 0 = (-1)^n\]
\\
We apply the operational polynomial $P(T)$ to each side of the equation:
:
\\
\begin{center}
$ P(T)e^{n \pi t i}=P(T)\Bigl(e^{n \pi i}\Bigr)^t=$ (use
\ref{EDF:Potencia})
\end{center}

\[=e^{n \pi t i}P(e^{n \pi i})=e^{n \pi t
i}P((-1)^n)=P((-1)^n)\cos n \pi t + i P((-1)^n) \sin n \pi t\]

and

\[P(T)e^{n \pi t i}=P(T)(\cos n \pi t + i \sin n \pi t)=P(T)\cos n \pi t + i P(T) \sin n \pi t\]
\\
Equating the real and imaginary parts and supposing that $P(-1) \ne 0$, we obtain

\begin{equation}
    P(T) \cos n \pi t=P(-1) \cos n \pi t \Rightarrow
    \mathbf{\frac{1}{P(T)}\cos n \pi t=\frac{\cos n \pi t}{P((-1)^n)}}
     \label{EDF:Coseno}
\end{equation}

\begin{equation}
P(T) \sin n \pi t=P(-1) \sin n \pi t \Rightarrow \mathbf{\frac{1}{P(T)}\sin n \pi
       t=\frac{\sin n \pi t}{P((-1)^n)}}
    \label{EDF:Seno}
\end{equation}

\begin{example}
Find a particular solution $y_P(t)$ of the finite difference equation
\[y(t+2)-5y(t+1)+6y(t)=\cos(\pi t)\]
\end{example}
\begin{solution}
\[y(t+2)-5y(t+1)+6y(t)=\cos(\pi t)\Rightarrow (T^2-5T+6)y(t)=\cos(\pi t)\]
Replacing \ref{EDF:Coseno}, with $n=1$, we obtain
\[ y(t)=\frac{1}{T^2-5T+6}\cos(\pi t)=\frac{\cos(\pi t)}{(-1)^2-5 \cdot (-1)+6} \Rightarrow y_P(t)=\frac{\cos(\pi t)}{12} \]
\end{solution}

\subsection{Equations where $\mathbf{\phi(t)=\lambda^t f(t)}$}

Consider the polynomial
\[P(T)=a_0 T^n+a_1T^{n-1}+\cdots+a_{n-1}T+a_n\]

We have
\[P(T)\lambda^t f(t)=a_0
T^n(\lambda^t f(t) )+a_1T^{n-1}(\lambda^t f(t) )
+\cdots+a_{n-1}T(\lambda^t f(t)) +a_n\lambda^t f(t)=\]
\[=a_0
\lambda^{t+n}f(t+n)+a_1\lambda^{t+n-1}f(t+n-1)+\cdots+a_{n-1}\lambda^{t+1}f(t+1)+a_n\lambda^tf(t)=\]
\[=\lambda^t(a_0
\lambda^{n}T^n
f(t)+a_1\lambda^{n-1}T^{n-1}f(t)+\cdots+a_{n-1}\lambda T f(t)+a_n
f(t)=\lambda^t P(\lambda T)f(t)\]
\\
Therefore, we obtain the property

\begin{equation}
   P(T)\lambda^t f(t)= \lambda^tP(\lambda T)f(t) \Leftrightarrow
   \mathbf{\frac{1}{P(T)}\lambda^t f(t)=\lambda^t\frac{1}{P(\lambda T)} f(t)} \ .
    \label{EDF:PotFun}
\end{equation}

\begin{example}
Find a particular solution $y_P(t)$ of the finite difference equation
\[y(t+2)-5y(t+1)+4y(t)=3^t \sin (\pi t) \ . \]
\end{example}
\begin{solution}
\[y(t+2)-5y(t+1)+4y(t)=3^t\sin (\pi t)\Rightarrow (T^2-5T+4)y(t)=3^t\sin (\pi t)\]
Replacing:
\[ y(t)=\frac{1}{T^2-5T+4}\Bigl[3^t\sin (\pi t)\Bigr]=3^t\frac{1}{(3T)^2-5(3T)+4}\sin (\pi t)=\]
\[ =3^t\frac{1}{9T^2-15T+4}\sin (\pi t)=3^t\frac{\sin (\pi t)}{9(-1)^2-15(-1)+4} \Rightarrow y_P(t)=3^t\frac{\sin (\pi t)}{28}\]
\end{solution}

\subsection{Formula for polynomials $\mathbf{P(T-\lambda)}$ with $\mathbf{\phi(t)=\lambda^t f(t)}$}

Here I prove by induction that $(T-\lambda)^n \lambda^t f(t)
= \lambda^t [\lambda(T-1)]^n f(t)$.

\begin{enumerate}

\item  When $n=1$
\[(T-\lambda) \lambda^t
f(t)=\lambda^{t+1}f(t+1)-\lambda^{t+1}f(t)=\lambda^{t+1}(T-1)f(t)=\lambda^{t}[\lambda(T-1)]f(t)\]

\item  Showing that the statement holds when $n=k$
\[(T-\lambda)^k \lambda^t f(t) = \lambda^t [\lambda(T-1)]^k f(t)\]

\item  Prove when $n=k+1$
\[(T-\lambda)^{k+1} \lambda^t
f(t)=(T-\lambda) (T-\lambda)^k \lambda^t f(t)  =(T-\lambda)
\lambda^{t}[\lambda(T-1)]^k f(t) =
\]
\[=\lambda^{t+1}[\lambda(T-1)]^k f(t+1)-\lambda^{t+1}[\lambda(T-1)]^k f(t)=\lambda^{t+1}[\lambda(T-1)]^{k} [Tf(t)-f(t)]=\]
\[=\lambda^{t+1} [\lambda(T-1)]^{k} (T-1) f(t) =\lambda^{t+1} \lambda^k(T-1)^{k+1}f(t)= \lambda^t [\lambda(T-1)]^{k+1} f(t)\]

Finally, we generalice for $P(T-\lambda) \lambda^t f(t)$.

Consider the polynomial $P(T)=a_0 T^n+a_1T^{n-1}+\cdots+a_{n-1}T+a_n$.

We have
\[P(T-\lambda)\lambda^t f(t)=\]
\[=a_0 (T-\lambda)^n(\lambda^t f(t) )+a_1(T-\lambda)^{n-1}(\lambda^t f(t))
+\cdots+a_{n-1}(T-\lambda)(\lambda^t f(t)) +a_n\lambda^t f(t)=\]
\[=a_0\lambda^{t}[r(T-1)]^nf(t)+a_1\lambda^{t}[r(T-1)]^{n-1}f(t)
+\cdots+a_{n-1}\lambda^{t}[r(T-1)]f(t)+a_n\lambda^t f(t)=\]
\[=\lambda^t P[\lambda (T-1)]f(t)\]
\\
Therefore, we can establish the property

\[P(T-r)\lambda^t f(t)=\lambda^tP[\lambda (T-1)]f(t) \Leftrightarrow\]
\begin{equation}
     \Leftrightarrow \mathbf{\frac{1}{P(T-\lambda)}\lambda^t f(t)=
       \lambda^t\frac{1}{P[\lambda (T-1)]} f(t)}
    \label{EDF:PolDes}
\end{equation}

\end{enumerate}

\begin{example}
Find a particular solution $y_P(t)$ of the finite difference equation
\[y(t+1)-2y(t)=2^t\]
\end{example}
\begin{solution}
\[y(t+1)-2y(t)=2^t\Rightarrow (T-2)y(t)=2^t  \Rightarrow y(t)=\frac{1}{T-2}2^t\]
We cannot apply \ref{EDF:Potencia}, since it would give a division
by zero. Instead, we write
\[ y(t)=\frac{1}{T-2}(2^t \cdot 1) \ .\]
By applying \ref{EDF:PolDes}, we obtain
\[ y(t)=2^t \frac{1}{2(T-1)}(1)=2^t \frac{1}{2}\frac{1}{(T-1)}(1)=2^{t-1} \frac{1}{T-1}(1)\]
Finally, using \ref{EDF:Unidad}:
\[y_P(t)=2^{t-1} t\]
\end{solution}

\section{Conclusions} 
This paper developed a technique for solving linear finite difference equations with
constant coefficients. In addition, it proves several fundamental properties 
(linearity, translation, unity and propagation) of the polynomial translation operator
and establishes formulae for solving several forms of FDE: $\lambda^t, \cos n\pi t, \sin n\pi t, 
\lambda^t$  and $\lambda^t f(t)$.


\end{document}